\documentclass[12pt, reqno]{amsart}

\usepackage{amssymb}  

\theoremstyle{plain}
\newtheorem{theorem}{Theorem}[section]
\newtheorem*{main}{Main Theorem}
\newtheorem{lemma}[theorem]{Lemma}
\newtheorem*{definition}{Definition}
\newtheorem*{proposition}{Proposition}
\newtheorem{corollary}[theorem]{Corollary}
\theoremstyle{remark}

\numberwithin{equation}{section}

\newcommand{\seclabel}[1]{\label{sec:#1}} 
\newcommand{\thmlabel}[1]{\label{thm:#1}} 
\newcommand{\lemlabel}[1]{\label{lem:#1}} 
\newcommand{\corlabel}[1]{\label{cor:#1}} 
\newcommand{\prplabel}[1]{\label{prp:#1}} 
\newcommand{\deflabel}[1]{\label{def:#1}} 

\newcommand{\secref}[1]{\ref{sec:#1}} 
\newcommand{\lemref}[1]{\ref{lem:#1}} 
\newcommand{\corref}[1]{\ref{cor:#1}} 

\newcommand{\Aut}{\mathcal{A}ut}		
\newcommand{\Ann}{\mathrm{Ann}}		
\newcommand{\End}{\mathcal{E}nd}	
\newcommand{\ZEnd}{\mathcal{ZE}nd}	
\newcommand{\QEnd}{\mathcal{QE}nd}	
\newcommand{\SEnd}{\mathcal{SE}nd}	

\newcommand{\Id}{\mathrm{Id}}	
\newcommand{\iv}{^{-1}}				
\title[F-Quasigroups and Generalized Modules]
{F-Quasigroups and Generalized Modules}

\author[T.~Kepka]{Tom\'{a}\v{s}~Kepka$^*$}
\thanks{$^*$Partially supported by the institutional grant
MSM 113200007 and by the Grant Agency of Charles University,
grant \#269/2001/B-MAT/MFF}
\address{Department of Algebra \\
MFF UK, Sokolovsk\'{a} 83 \\
186 75 Praha 8, Czech Republic}
\email{kepka@karlin.mff.cuni.cz}
\author[M.~K.~Kinyon]{Michael~K.~Kinyon}
\address{Department of Mathematical Sciences \\
Indiana University South Bend \\
South Bend, IN 46634 USA}
\email{mkinyon@iusb.edu}
\urladdr{http://mypage.iusb.edu/\symbol{126}mkinyon}
\author[J.~D.~Phillips]{J.~D.~Phillips}
\address{Department of Mathematics \& Computer Science \\
Wabash College \\
Crawfordsville, IN 47933 USA}
\email{phillipj@wabash.edu}
\urladdr{http://www.wabash.edu/depart/math/faculty.html{\#}Phillips}

\date{\today}


\begin{document}

\begin{abstract}
In \cite{KKP}, we showed that every F-quasigroup is linear over a
special kind of Moufang loop called an NK-loop. Here we extend this
relationship by showing an equivalence between the equational class
of (pointed) F-quasigroups and the equational class corresponding to a  
certain notion of generalized module (with noncommutative, nonassociative
addition) for an associative ring.
\end{abstract}

\maketitle

\section{Introduction}
\seclabel{intro}

A \emph{quasigroup} $(Q,\cdot)$ is a set $Q$ with a binary operation
$\cdot : Q\times Q\to Q$, denoted by juxtaposition, such that for
each $a,b\in Q$, the equations $ax = b$ and $ya = b$ have unique
solutions $x,y\in Q$. 
In a quasigroup $(Q,\cdot)$, there exist transformations
$\alpha, \beta : Q\to Q$ such that $x \alpha(x) = x = \beta(x) x$ for all $x\in Q$.
Now $(Q,\cdot)$ is called an \emph{F-quasigroup} if it satisfies the equations
\[
x \cdot yz = y \cdot \alpha (x) z
\qquad \text{and} \qquad
zy \cdot x = z \beta (x) \cdot yx
\]
for all $x,y,z \in Q$.

If $(Q,\cdot)$ is a quasigroup, we set
$M(Q) = \{a \in Q : xa \cdot yx = xy \cdot ax,\ \forall x,y \in Q \}$.
If $(Q,\cdot)$ is an F-quasigroup, then $M(Q)$ is a normal
subquasigroup of $Q$ and $Q/M(Q)$ is a group \cite[Lemma 7.5]{KKP}.

We denote by $\mathcal{F}_p$ the equational class (and category)
of pointed $F$-quasigroups. That is, $\mathcal{F}_p$ consists of ordered pairs
$(Q,a)$, where $Q$ is an F-quasigroup and $a \in Q$. We put
$\mathcal{F}_m = \{(Q,a) \in \mathcal{F}_p : a \in M(Q)\}$. 

A quasigroup with a neutral element is called a \emph{loop}. Throughout
this paper, we adopt an additive notation convention $(Q,+)$ (with
neutral $0$) for loops,
although we do not assume that $+$ is commutative. The \emph{nucleus}
of a loop $(Q,+)$ is the set
\[
N(Q,+) = \{ a\in Q :
\left\{\begin{array}{c}
(a + x) + y = a + (x + y) \\
(x + a) + y = x + (a + y) \\ 
(x + y) + a = x + (y + a)
\end{array}
\right\},
\forall x,y\in Q \} .
\]
The \emph{Moufang center} is the set
\[
K(Q,+) = \{ a\in Q : (a + a) + (x + y) = (a + x) + (a + y), \; \ \forall x,y\in Q \} .
\]
The intersection of the nucleus and Moufang center of a loop is the
\emph{center} $Z(Q,+) = N(Q,+)\cap K(Q,+)$.
Each of the nucleus, the Moufang center, and the center is a
subloop, and the center is, in fact, a normal subloop \cite{Br,Pf}.

A $(Q,+)$ will be called an \emph{NK-loop} if for each $x \in Q$,
there exist $u \in N(Q,+)$ and $v \in K(Q,+)$ such that
$x = u + v \; (= v + u)$. In other words, $Q$ can be decomposed
as a central product $Q = N(Q,+)\; K(Q,+)$. It was shown in
\cite{KKP} that every NK-loop is a Moufang A-loop. A \emph{Moufang loop} is
a loop satisfying the identity $((x + y) + x) + z = x + (y + (x + z))$
or any of its known equivalents \cite{Br,Pf}. Every Moufang loop is
\emph{diassociative}, that is, the subloop generated by any given
pair of elements is a group \cite{moufang}. For a loop $(Q,+)$, the
\emph{inner mapping group} is the stabilizer of $0$ in the group
of permutations of $Q$ generated by all left and right translations
$L_x y = x + y = R_y x$. An \emph{A-loop} is a loop such
that every inner mapping is an automorphism \cite{BP}.

In any Moufang A-loop $(Q,+)$, such as an NK-loop, the nucleus
$N(Q,+)$ is normal (in fact, this is true in any Moufang
loop), and $Q/N(Q,+)$ is a commutative Moufang loop of exponent $3$.
In particular, for each $x\in Q$, $3x\in N(Q,+)$, where $3x = x + x + x$.
The Moufang center $K(Q,+)$ is also normal in $Q$ (but this is not
necessarily the case in arbitrary Moufang loops), and 
$Q/K(Q,+)$ is a group \cite[Lemma 4.3]{KKP}. In an NK-loop $(Q,+)$, we also have
$Z(Q,+) = Z(N(Q,+)) = K(N(Q,+)) = Z(K(Q,+)) = N(K(Q,+))$. 
In addition, $K(Q,+) = \{ a \in Q : a + x = x + a \; \ \forall x\in Q \}$.

The connection between F-quasigroups and NK-loops was established in
\cite{KKP}.

\begin{proposition}
\prplabel{kkp}
For a quasigroup $(Q,\cdot)$, the following are equivalent:
\begin{enumerate}
\item[1.] $(Q,\cdot)$ is an F-quasigroup.
\item[2.] There exist an NK-loop $(Q,+)$,  
$f,g\in \Aut(Q,+)$, and $e\in N(Q,+)$ such that
$x \cdot y = f(x) + e + g(y)$
for all $x,y\in Q$, $fg = gf$, and
$x + f(x), x + g(x)\in N(Q,+)$,
$-x + f(x),-x + g(x)\in K(Q,+)$ for all $x\in Q$.
\end{enumerate}
\end{proposition}

We refer to the data $(Q,+,f,g,e)$ of the proposition as being
an \emph{arithmetic form} of the F-quasigroup $(Q,\cdot)$. 
If $(Q,a)$ is a pointed F-quasigroup in $\mathcal{F}_p$, then there
is an arithmetic form such that $a = 0$ is the neutral element of $(Q,+)$.

The purpose of this paper is to extend the connection between (pointed) F-quasigroups
and NK-loops further by showing an equivalence of equational classes between 
$\mathcal{F}_p$ and a certain notion of generalized module for an associative ring.
Thus the study of (pointed) F-quasigroups effectively becomes a part of ring theory. 
The generalization we require weakens the additive abelian group structure of
a module to an NK-loop structure.

\begin{definition}
\deflabel{gen-module}
Let $R$ be an associative ring, possibly without unity. A
\emph{generalized (left)} $R$\emph{-module} is an $NK$-loop $(Q,+)$ supplied with an
$R$-scalar multiplication $R \times Q \rightarrow Q$ such that the following conditions
are satisfied: for all $a,b \in R$, $x, y \in Q$, $z \in N(Q,+)$, and $w \in K(Q,+)$,
\begin{enumerate}
\item[1.] \qquad $a(x + y) = ax + ay$,
\item[2.] \qquad $(a + b)x = ax + bx$,
\item[3.] \qquad $a(bx) = (ab)x$,
\item[4.] \qquad $ax \in K(Q,+)$,
\item[5.] \qquad $az \in N(Q,+)$, and
\item[6.] \qquad there exists an integer $m$ such that $mw + aw \in Z(Q,+)$.
\end{enumerate}
\end{definition}

Here $mw = w + \cdots + w$ ($m$ terms) is unambiguous by diassociativity.

If $Q$ is a generalized $R$-module, then define the \emph{annihilator} of
$Q$ to be $\Ann(Q) = \{a \in R : aQ = 0\}$. Clearly, $\Ann(Q)$ is an ideal
of the ring $R$.

In order to state our main result, we need to describe a particular ring.
Let $\textbf{S} = \mathbb{Z}[\textbf{x},\textbf{y},\textbf{u},\textbf{v}]$
be the polynomial ring in four commuting indeterminates 
$\textbf{x},\textbf{y},\textbf{u}$, and $\textbf{v}$ over the ring 
$\mathbb{Z}$ of integers. Put
$\textbf{R} = \textbf{S} \textbf{x} + \textbf{S} \textbf{y} + \textbf{S} \textbf{u} + \textbf{S} \textbf{v}$,
so that $\textbf{R}$ is the ideal generated by the indeterminates. Clearly, $\textbf{R}$
is a free commutative and associative ring (without unity) freely generated by the indeterminates.

Let $\mathcal{M}$ denote the equational class (and category) of generalized
$\textbf{R}$-modules $Q$ such that:
\begin{enumerate}
\item[1.] \qquad $2z + \textbf{x}z \in N(Q,+)$, $2z + \textbf{y}z \in N(Q,+)$ for all $z \in Q$,
\item[2.] \qquad $\textbf{x} + \textbf{u} + \textbf{x} \textbf{u} \in \mathrm{Ann}(Q)$, and
\item[3.] \qquad $\textbf{y} + \textbf{v} + \textbf{y} \textbf{v} \in \mathrm{Ann}(Q)$.
\end{enumerate}
Further, let $\mathcal{M}_p$ be the equational class (and category) of pointed 
objects from $\mathcal{M}$. That is, $\mathcal{M}_p$ consists of ordered pairs $(Q,e)$,
where $Q \in \mathcal{M}$ and $e \in Q$.
Put $\mathcal{M}_n = \{ (Q,e)\in \mathcal{F}_p : e\in N(Q,+)\}$, the equational
class (and category) of \emph{nuclearly} pointed objects from $\mathcal{M}$, and put
$\mathcal{M}_c = \{ (Q,e)\in \mathcal{F}_p : e\in Z(Q,+)\}$, the equational
class (and category) of \emph{centrally} pointed objects from $\mathcal{M}$.

Our main result is the following equivalence between pointed F-quasigroups
and generalized $\textbf{R}$-modules.

\begin{main}
\thmlabel{main}
The equational classes $\mathcal{F}_p$ and $\mathcal{M}_n$ are equivalent.
The equivalence restricts to an equivalence between $\mathcal{F}_m$ and $\mathcal{M}_c$
\end{main}

\section{Quasicentral endomorphisms}
\seclabel{quasi-end}

In this section, let $(Q,+)$ denote a (possibly non-commutative) diassociative loop.
We endow the set $\End(Q,+)$ of all endomorphisms of $(Q,+)$ with the standard
operations of addition, negation, and composition, \textit{viz.},
for $f,g\in \End(Q,+)$, $f + g$ is defined by $(f + g)(x) = f(x) + g(x)$, $-f$ is
defined by $(-f)(x) = -f(x) = f(-x)$, and $fg$ is defined by $fg(x) = f(g(x))$
for all $x\in Q$.

An endomorphism $f$ of $(Q,+)$ is called \emph{central} if $f(Q) \subset Z(Q,+)$.
We denote the set of all central endomorphisms of $(Q,+)$ by $\ZEnd(Q,+)$.
The verification of the following result is easy and omitted.

\begin{lemma}
\lemlabel{cend-ring}
Let $f, g, h\in \ZEnd(Q,+)$ be given. Then:
\begin{enumerate}
\item[1.] $f + g\in \ZEnd(Q,+)$,
\item[2.] $f + (g + h) = (f + g) + h$ and $f + g = g + f$,
\item[3.] the zero endomorphism of $(Q,+)$ is central,
\item[4.] $-f\in \ZEnd(Q,+)$,
\item[5.] $f + (-f) = 0$ and $f + 0 = f$.
\item[6.] $fg\in \ZEnd(Q,+)$,
\end{enumerate}
\end{lemma}

\begin{corollary}
\corlabel{cend-ring}
$\ZEnd(Q,+)$ is an associative ring (possibly without unity)
with respect to the standard operations.
\end{corollary}

Let $m$ be an integer. An endomorphism $f$ of $(Q,+)$ is called
$m$-\emph{quasicentral} if $mx + f(x) \in Z(Q,+)$ for all
$x \in Q$ (in which case $mx + f(x) = f(x) + mx$).
An endomorphism is called \emph{quasicentral} if it is $m$-quasicentral
for at least one integer $m$. We denote by $\QEnd(Q,+)$
the set of all quasicentral endomorphisms of $(Q,+)$.
The following is an obvious consequence of these definitions.

\begin{lemma}
\lemlabel{quasi}
\begin{enumerate}
\item[1.] An endomorphism is $0$-quasicentral if and only if it is central,
\item[2.] $\ZEnd(Q,+) \subset \QEnd(Q,+)$, and
\item[3.] the identity automorphism, $\Id_Q$, of $(Q,+)$ is $(-1)$-quasicentral.
\end{enumerate}
\end{lemma}

\begin{lemma}
\lemlabel{quasi2}
Let $f,g\in \End(Q,+)$.
\begin{enumerate}
\item[1.] If $f$ is $m$-quasicentral and $g$ is $n$-quasicentral, then
$fg$ is ($-mn$)-quasicentral.
\item[2.] If $f,g\in \QEnd(Q,+)$, then $fg\in \QEnd(Q,+)$.
\end{enumerate}
\end{lemma}

\begin{proof}
For (1): Fix $x\in Q$.
Since $f$ is $m$-quasicentral, $g(mx) + fg(x) = mg(x) + fg(x) \in Z(Q,+)$.
Since $g$ is $n$-quasicentral, $-mnx -mg(x) = -(g(mx) + nmx) \in Z(Q,+)$.
Consequently,
\begin{align*}
-mnx + fg(x) &= ([-mnx - mg(x)] + mg(x)) + fg(x) \\
&= [-mnx - mg(x)] + [mg(x) + fg(x)] \in Z(Q,+).
\end{align*}
Thus, $fg$ is ($-mn$)-quasicentral, as claimed.

(2) follows immediately from (1).
\end{proof}

\begin{lemma}
\lemlabel{quasi-comm}
Assume that $(Q,+)$ is commutative, let $f,g\in \End(Q,+)$ be $m$-quasicentral
and $n$-quasicentral, respectively. Then
\begin{enumerate}
\item[1.] \qquad $-f$ is $(-m)$-quasicentral,
\item[2.] \qquad $f + g$ is an $(m + n)$-quasicentral endomorphism.
\end{enumerate}
In particular, for $f,g\in \QEnd(Q,+)$, $-f, f+g\in \QEnd(Q,+)$.
\end{lemma}

\begin{proof}
(1) is clear. For (2), 
set $z = (-mx - f(x)) + (-my - f(y)) + (-mx - g(x)) + (-my - g(y))$.
Then $z\in Z(Q,+)$. It follows that
\begin{align*}
z + (f + g)(x + y) &= z + ([f(x) + f(y)] + [g(x) + g(y)]) 
= -mx - my - nx - ny \\
&= -mx - nx - my - ny = z + ([f(x) + g(x)] + [f(y) + g(y)])\\
&= z + ((f + g)(x) + (f + g)(y)) .
\end{align*}
Thus $f + g\in \End(Q,+)$. Similarly, $(m + n)x + (f + g)(x) = [mx + f(x)]
+ [nx + g(x)] \in Z(Q,+)$. That is, (2) holds.
\end{proof}

\begin{lemma}
\lemlabel{quasi-comm2}
Assume that $(Q,+)$ is commutative and let $f,g,h\in \QEnd(Q,+)$. Then
\begin{enumerate}
\item[1.] \qquad $f + g = g + f$,
\item[2.] \qquad $f + (g + h) = (f + g) + h$,
\item[3.] \qquad $f + (-f) = 0$, and
\item[4.] \qquad $f + 0 = f$.
\end{enumerate}
\end{lemma}

\begin{proof}
(1), (3), and (4) are obvious. For (2): 
There exist $m, n, k\in \mathbb{Z}$ such that $mx + f(x), nx + g(x),
kx + h(x)\in Z(Q,+)$. Set
$y = (-f(x) - mx) + (-g(x) - nx) + (-h(x) - kx)$. Then $y \in Z(Q,+)$
and $y + (f(x) + (g(x) + h(x)))
= -(m + n + k)x = y + ((f(x) + g(x)) + h(x))$ for all $x \in Q$.
\end{proof}

\begin{corollary}
\corlabel{quasi-comm}
If $(Q,+)$ is commutative, then $\QEnd(Q,+)$ is an associative ring with unity. 
\end{corollary}

We conclude this section with a straightforward observation.

\begin{lemma}
\lemlabel{quasi-3}
Assume that for $k\in \{ 1,2,3\}$, $kx \in Z(Q,+)$ for all $x \in Q$. Then
\begin{enumerate}
\item[1.] Every quasicentral endomorphism is $m$-central for some
$m \in \{0,1,-1\}$,
\item[2.] If $f\in \QEnd(Q,+)\cap \Aut(Q,+)$, then $f\iv \in \QEnd(Q,+)$.
\end{enumerate}
\end{lemma}

\section{Special endomorphisms of $NK$-loops}
\seclabel{special}

In this section, let $(Q,+)$ be an $NK$-loop. We denote
by $N$, $K$, and $Z$ the underlying sets of
$N(Q,+)$, $K(Q,+)$, and $Z(Q,+)$, respectively. As noted in
{\S}\secref{intro}, $Z(Q,+) = Z(N,+) = Z(K,+)$ and $Z = N \cap K$.

An endomorphism $f$ of $(Q,+)$ will be called \emph{special} if 
$f(Q) \subset K$, $f|_K$ is a quasicentral endomorphism of $(K,+)$,
and $f(N) \subset N$. Then $f|_N$ is a central endomorphism of $(N,+)$
and $f(N) \subset Z$. 
We denote by $\SEnd(Q,+)$ the set of special endomorphisms of $(Q,+)$.

\begin{lemma}
\lemlabel{special}
Let $f,g,h\in \SEnd(Q,+)$. Then
\begin{enumerate}
\item[1.] \qquad $fg\in \SEnd(Q,+)$,
\item[2.] \qquad $f+g\in \SEnd(Q,+)$, and $f + g = g + f$,
\item[3.] \qquad $f + (g + h) = (f + g) + h$,
\item[4.] \qquad $-f\in \SEnd(Q,+)$, $f + (-f) = 0$, and $f + 0 = f$.
\end{enumerate}
\end{lemma}

\begin{proof}
For (1), use Lemma \lemref{quasi2}. 

For (2):
Take $x,y \in Q$. Then $x = a + b$ and $y = c + d$ for some
$a,c \in N, b,d \in K$ so that
\begin{align*}
u &= (f + g)(x + y) = f(x + y) + g(x + y) = [f(x) + f(y)] + [g(x) + g(y)] \\
&= [(f(a) + f(b)) + (f(c) + f(d))] + [(g(a) + g(b)) + (g(c) + g(d))]
\end{align*}
and
\begin{align*}
v &= (f + g)(x) + (f + g)(y) = [f(x) + g(x)] + [f(y) + g(y)] \\
&= [(f(a) + f(b)) + (g(a) + g(b))] + [(f(c) + f(d)) + (g(c) + g(d))].
\end{align*}
The restrictions $f|_N$ and $g|_N$ are central endomorphisms of $(N,+)$, and it follows that
$f(N) \cup g(N) \subset Z(N,+) = Z(Q,+)$. Thus, $f(a), f(c), g(a), g(c) \in Z$ and
in order to check that $u = v$ it is sufficient to show that $(f(b) + f(d)) + (g(b) + g(d)) 
= (f(b) + g(b)) + (f(d) + g(d))$. However, the latter equality holds, since the
restrictions $f|_K$ and $g|_K$ are quasicentral endomorphisms of the commutative loop
$(K,+)$ and Corollary \corref{quasi-comm} applies.

We have shown that $f + g\in \End(Q,+)$. The facts that $f + g$ is special
and $f + g = g + f$ are easily seen, using Lemma \lemref{quasi-comm2} applied to the
loop $(K,+)$.

For (3): 
Using the facts that $(Q,+)$ is an NK-loop and
$f(N) \cup g(N) \cup h(N) \subset Z$, it is enough to show that
$f(u) + (g(u) + h(u)) = (f(u) + g(u)) + h(u)$
for all $u \in K$. Now we proceed similarly as in the proof of Lemma \lemref{quasi-comm2}.

Finally, (4) is easy.
\end{proof}

\begin{corollary}
\corlabel{send-ring}
$\SEnd(Q,+)$ is an associative ring (possibly without unity).
\end{corollary}

An endomorphism $f$ of $(Q,+)$ will be said to satisfy \emph{condition (F)} if
\[ 
-x + f(x) \in K \qquad \text{and} \qquad x + f(x) \in N
\]
for all $x \in Q$. Then $f(K) \subset K$ and $f(N) \subset N$.

\begin{lemma}
\lemlabel{h-send}
Let $f\in \End(Q,+)$ satisfy \emph{(F)}. Define 
$h : Q\to Q$ by $h(x) = -x + f(x)$ for all $x \in Q$. Then $h\in \SEnd(Q,+)$.
\end{lemma}

\begin{proof}
First we check that $h\in \End(Q,+)$.
Fix $x,y\in Q$ with $x = a + b$, $y = c + d$,
$a,c \in N$, $b,d \in K$.
Set $u = h(x+y)$, $v = h(x)+h(y)$, and 
$w = (a - f(a)) + (c - f(c)) + (-b - f(b)) + (-d - f(d))$.
Then $w\in Z$,
\[
u = (-y - x) + f(x + y) = ((-d - c) + (-b - a))
+ ((f(a) + f(b))) + (f(c) + f(d))
\]
and
\[
v = (-x + f(x)) + 
(-y + f(y)) = (-b - a) + (f(a) + f(b)) + ((-d - c) + (f(c) + f(d))) .
\]
On the other hand,
\begin{align*}
u + w &= [(-d - c) + (-b - a)] + [(a - b) + (c - d)] \\
&= [(-d - c) - b] + [-a + (a - b)] + (c - d)  \\
&= [(-d - c) - b] + [(c - d) - b] \\
&= [(-d - c) + (c - d)] - 2b \\
&= -2d - 2b \\
&= -2(b + d) \\
&= [(-b - a) + (a - b)] + [(-d - c) + (c - d)] \\
&= v + w. 
\end{align*}
Consequently, $u = v$, so that $h\in \End(Q,+)$, as claimed.
Further, it follows immediately from the definition
of $h$ that $h(Q) \subset K$ and $h(N) \subset N$ (then $h(N) \subset Z$).
Finally, $2a + h(a) = a + f(a) \in Z$ for all $a \in K$, and therefore
$h|_K$ is a $2$-quasicentral endomorphism of $(K,+)$. Thus $h\in \SEnd(Q,+)$.
\end{proof}

\begin{lemma}
\lemlabel{send-commute}
Let $f, g\in \End(Q,+)$ satisfy \emph{(F)}. 
Define $h,k : Q\to Q$ by $h(x) = -x + f(x)$ and $k(x) = -x + g(x)$
for all $x \in Q$. Then $hk = kh$ if and only if $fg = gf$.
\end{lemma}

\begin{proof}
By Lemma \lemref{h-send}, $h\in \End(Q,+)$, and hence
\[
hk(x) = h(-x + g(x)) = -h(x) + hg(x) = (-f(x) + x) + (-g(x) + fg(x)).
\]
On the other hand, 
\[
kh(x) = -h(x) + gh(x) = (-f(x) + x) + (-g(x) + gf(x))
\]
by the definition of $h$ and $k$. The result is now clear.
\end{proof}

\begin{lemma}
\lemlabel{aut-F}
Let $f,g\in \Aut(Q,+)$ satisfy \emph{(F)}. Define $h,k,p,q: Q\to Q$
by $h(x) = -x + f(x)$, $k(x) = -x + g(x)$, $p(x) = -x + f^{-1}(x)$, and
$q(x) = -x + g^{-1}(x)$ for all $x \in Q$. Then
\begin{enumerate}
\item[1.] $h,k,p,q\in \SEnd(Q,+)$,
\item[2.] $hp = ph$ and $h + p + hp = 0$,
\item[3.] $kq = qk$ and $k + q + kq = 0$, and
\item[4.] If $fg = gf$, then the endomorphisms $h,k,p,q$ commute pairwise.
\end{enumerate}
\end{lemma}

\begin{proof}
(1) follows from Lemma \lemref{h-send}.

For (2): 
We have $ff\iv = f\iv f$ and hence $hp = ph$ by Lemma \lemref{send-commute}.
Now, put $A = h + p + hp$. Then $A$ is a (special) endomorphism of $(Q,+)$ and
$A(x) = [-x + f(x)] + [-x + f^{-1}(x)] + [(-f^{-1}(x) + x) + (-f(x) + x)]$.
Clearly, $N \subset \ker(A) ( = \{u \in Q : A(u) = 0\} )$.
On the other hand, if $x \in K$, then $-x + f(x), -x + f^{-1}(x) \in Z$ and the equality
$A(x) = 0$ is clear, too. Thus, $N \cup K \subset \ker(A)$. But $(Q,+)$ is an $NK$-loop and 
$\ker(A)$ is a subloop of $(Q,+)$. It follows $\ker(A) = Q$ and $A = 0$.

(3) is proven similarly to (2).

For (4), combine (2), (3), and Lemma \lemref{send-commute}.
\end{proof}

\section{The equivalence}
\seclabel{equiv}

We now turn to the proof of the Main Theorem. First, recall the definition
of generalized module over a ring $R$, and observe that the conditions
(1), (4), (5), and (6) of the definition imply that for each $a\in R$,
the transformation $Q\to Q; x \mapsto ax$ is a special endomorphism
of $(Q,+)$. Recall also the ring $\textbf{R}$, which is the ideal of
$\textbf{S} = \mathbb{Z}[\textbf{x},\textbf{y},\textbf{u},\textbf{v}]$
freely generated by the commuting indeterminates 
$\textbf{x},\textbf{y},\textbf{u}$, and $\textbf{v}$.

First, take $(Q,a) \in \mathcal{F}_p$. As described in {\S}\secref{intro},
let $(Q,+,f,g,e)$ be the arithmetic form of the F-quasigroup $(Q,\cdot)$
such that $a = 0$ in $(Q,+)$. Then $f,g\in \Aut(Q,+)$ satisfy condition \emph{(F)}.
Further, define $\varphi, \mu, \psi, \nu : Q\to Q$
by $\varphi(x) = -x + f(x)$, $\mu(x) = -x + f^{-1}(x)$, $\psi(x) = -x + g(x)$, and
$\nu(x) = -x + g^{-1}(x)$ for all $x \in Q$.
By Lemma \lemref{aut-F}, the special endomorphisms $\varphi, \psi, \mu$, and $\nu$
of the loop $(Q,+)$ commute pairwise, and
$\varphi + \mu + \varphi \mu = 0 = \psi + \nu + \psi \nu$. Consequently,
these endomorphisms generate a commutative subring of the ring $\SEnd(Q,+)$ 
(see Corollary \corref{send-ring}) and there exists a (uniquely determined) homomorphism
$\lambda : \textbf{R} \to \SEnd(Q,+)$ such that $\lambda(\textbf{x}) = \varphi$, 
$\lambda(\textbf{y}) = \psi$, $\lambda(\textbf{u}) = \mu$, and 
$\lambda(\textbf{v}) = \nu$. The homomorphism $\lambda$ induces an $\textbf{R}$-scalar
multiplication on the loop $(Q,+)$, and the resulting generalized $\textbf{R}$-module 
will be denoted by $\overline{Q}$. By Lemma \lemref{aut-F}, $\lambda(\textbf{x} + \textbf{u} +
\textbf{x} \textbf{u}) = 0 = \lambda(\textbf{y} + \textbf{v} + \textbf{y} \textbf{v})$,
and so
$\textbf{x} + \textbf{u} + \textbf{x} \textbf{u},
\textbf{y} + \textbf{v} + \textbf{y} \textbf{v} \in \Ann(Q)$.
Also, since $f,g$ satisfy \emph{(F)}, we have $2z + \lambda(\textbf{x})z = 2z + \varphi(z)
= z + f(z)\in N(Q,+)$ and similarly $2z + \lambda(\textbf{y})z\in N(Q,+)$ for all $z\in Q$.
It follows that $\overline{Q} \in \mathcal{M}$. Now define
$\rho : \mathcal{F}_p \to \mathcal{M}_n$ by 
$\rho(Q,a) = (\overline{Q},e)$, and observe that 
$(\overline{Q},e) \in \mathcal{M}_c$ if and only if $e \in Z(Q,+)$.

Next, take $(\overline{Q},e) \in \mathcal{M}_n$ and define $f,g : Q\to Q$ by
$f(z) = z + \textbf{x} z$ and $g(z) = z + \textbf{y} z$ for all $z \in Q$.
We have $f(x + y) = (x + y) + (\textbf{x} x + \textbf{x} y)$ and
$f(x) + f(y) = (x + \textbf{x} x) + (y + \textbf{x} y)$. Further,
$x = u_1 + v_1$, $y = u_2 + v_2$ for some $u_1, u_2 \in N(Q,+)$,
$v_1, v_2 \in K(Q,+)$, and hence,
$f(x + y) = (u_1 + u_2 + v_1 + v_2) + (\textbf{x} u_1 + \textbf{x} u_2 + \textbf{x} v_1 + \textbf{x} v_2)$,
and
$f(x) + f(y) =
(u_1 + \textbf{x} u_1 + v_1 + \textbf{x} v_1) + (u_2 + \textbf{x} u_2 + v_1 + \textbf{x} v_2)$.
But
$\textbf{x} u_1, \textbf{x} u_2 \in Z(Q,+)$, and so in order to show
$f(x + y) = f(x) + f(y)$, it is enough to check that
$(v_1 + v_2) + (\textbf{x} v_1 + \textbf{x} v_2) =
(v_1 + \textbf{x} v_1) + (v_2 + \textbf{x} v_2)$. However,
$-2v_1 - \textbf{x} v_1 \in Z(Q,+)$ and 
$-2v_2 - \textbf{x} v_2 \in Z(Q,+)$, and so the latter equality is clear.

We have proven that $f\in \End(Q,+)$, and the proof that $g\in \End(Q,+)$
is similar. Now by definition of generalized module,
$-x + f(x) = \textbf{x} x \in K(Q,+)$
and $-x + g(x) = \textbf{y} x \in K(Q,+)$ for all $x \in Q$.
By definition of $\mathcal{M}$, $x + f(x) = 2x + \textbf{x} x \in N(Q,+)$ and
$x + g(x) = 2x + \textbf{y} x \in N(Q,+)$ for all $x \in Q$.
This means that both $f$ and $g$ satisfy \emph{(F)} and it follows
from Lemma \lemref{send-commute} that $fg = gf$.

Define $h: Q\to Q$ by $h(x) = x + \textbf{u} x$ for $x\in Q$.
We have $\textbf{u} x + \textbf{x} x + \textbf{x} \textbf{u} x = 0$,
and so  $\textbf{x} x + \textbf{x} \textbf{u} x = -\textbf{u} x$.
Now, $fh(x) = h(x) + \textbf{x} h(x) = (x + \textbf{u} x) +
(\textbf{x} x + \textbf{x} \textbf{u} x) = (x + \textbf{u} x) - \textbf{u} x
= x$ and $fh = \Id_Q$. Similarly, $hf = \Id_Q$ and we see that $f\in \Aut(Q,+)$.
Similarly, $g\in \Aut(Q,+)$.

We have that $f, g\in \Aut(Q,+)$, and $e\in Q$ satisfy the conditions
of the Proposition of {\S}\secref{intro}, and so defining a multiplication on $Q$ by
$xy = f(x) + e + g(y)$ for all $x,y \in Q$ gives an $F$-quasigroup.
Define $\sigma : \mathcal{M}_n \to \mathcal{F}_p$ by 
$\sigma(\overline{Q},e) = (Q,0)$.

It is easy to check that the operators $\rho$ and $\sigma$
represent an equivalence between $\mathcal{F}_p$ and $\mathcal{M}_n$.
Further, $0 \in \mathcal{M}(Q)$ if and only if $e \in Z(Q,+)$, so
that $\rho$ and $\sigma$ restrict to an equivalence between
$\mathcal{F}_m$ and $\mathcal{M}_c$. This completes the proof
of the Main Theorem.

\end{document}